\documentclass[12pt,reqno]{article}
\usepackage{amssymb}
\usepackage{amsmath}
\usepackage{amsthm}
\usepackage{amssymb}
\usepackage{amsmath}
\usepackage{amsthm}
\newcommand{\br}[3]{{$#1$}$\lower4pt\hbox{$\tp\atop\raise4pt \hbox{$\scriptscriptstyle{#2}$}$} ${$#3$}}
\newcommand{\tw}[3]{{$#1$}${\,\scriptscriptstyle {#2}}\atop\raise9pt\hbox{$\scriptstyle\tp$} ${$#3$}}
\newcommand{\ttps}[2]{{#1}\raise5pt\hbox{$\lower12pt\hbox{$\scriptstyle\tp$}\atop \lower0pt\hbox{$\tilde\;$}$}\raise4.5pt\hbox{${\scriptstyle{#2}}$}}

\newcommand{\ttp}{{\lower12pt\hbox{$\tp$}\atop \hbox{$\tilde\;$}}}

\newcommand{\btr}{\raise1.2pt\hbox{$\scriptstyle\blacktriangleright$}\hspace{2pt}}

\newcommand{\tp}{\otimes}
\renewcommand{\L}{\mathcal{L}}
\newcommand{\A}{\mathcal{A}}

\newcommand{\J}{\mathcal{J}}
\newcommand{\N}{\mathbb{N}}
\newcommand{\C}{\mathbb{C}}
\newcommand{\R}{\mathbb{R}}
\newcommand{\Z}{\mathbb{Z}}

\newcommand{\Ha}{\mathcal{H}}

\newcommand{\Rg}{\mathfrak{R}}

\newcommand{\V}{V}
\newcommand{\U}{\mathcal{U}}

\newcommand{\gm}{\gamma}

\newcommand{\la}{\lambda}
\newcommand{\End}{\mathrm{End}}
\newcommand{\Tr}{\mathrm{Tr}}
\newcommand{\tr}{\triangleright}

\newcommand{\g}{\mathfrak{g}}
\renewcommand{\b}{\mathfrak{b}}
\newcommand{\h}{\mathfrak{h}}
\newcommand{\n}{\mathfrak{n}}
\newcommand{\nn}{\nonumber}
\newcommand{\mb}{{\bf m}}
\newcommand{\p}{\mathfrak{p}}
\renewcommand{\l}{\mathfrak{l}}
\renewcommand{\c}{\mathfrak{c}}

\newcommand{\si}{\sigma}
\newcommand{\al}{\alpha}
\newcommand{\bt}{\beta}
\newcommand{\be}{\begin{eqnarray}}
\newcommand{\ee}{\end{eqnarray}}

\newtheorem{thm}{Theorem}[section]
\newtheorem{propn}[thm]{Proposition}
\newtheorem{lemma}[thm]{Lemma}

\newtheorem{remark}[thm]{Remark}
\newtheorem{definition}[thm]{Definition}
\newtheorem{example}[thm]{Example}

\theoremstyle{definition}
\newtheorem{parag}{}[section]

\begin{document}
\title{Quantum coadjoint orbits of $GL(n)$ and generalized Verma modules}
\author{J. Donin\footnote{Supported in part
by the Israel Academy of Sciences grant no. 8007/99-03.}\hspace{5pt}
 and A. I. Mudrov\footnote{Supported in part
by the Israel Academy of Sciences grant no. 8007/99-03 and by the RFBR grant no. 02-01-00085.}}
\setcounter{page}{1}
\date{}
\maketitle
\vspace{-0.0 in}
\begin{center}
{\small Department of Mathematics, Bar Ilan
University, 52900 Ramat Gan, Israel.}
\end{center}
\begin{abstract}
In our previous paper, we constructed an explicit $GL(n)$-equivariant
quantization of the Kirillov--Kostant-Souriau bracket on a semisimple coadjoint orbit.
In the present paper, we realize that quantization as a subalgebra of endomorphisms of
a generalized Verma module. As a corollary, we obtain an explicit description of the annihilators
of generalized Verma modules over $\U\bigl(gl(n)\bigr)$. As an application, we construct real forms
of the quantum orbits and classify finite dimensional representations.
We compute the non-commutative Connes index for basic homogenous vector bundles over the quantum orbits.
\\[16pt]
{\small \underline{Key words}: Kirillov-Kostant-Souriau bracket, equivariant quantization, generalized Verma modules.\\
\underline{AMS classification codes}: 53D55, 53D05, 22E47}
\end{abstract}
\section{Introduction}
Amongst Poisson manifolds with symmetries the simplest and the most natural one is the dual space $\g^*$ to
a Lie algebra $\g$.  The Poisson structure on $\g^*$ is induced by the Lie bracket of $\g$ considered
as a bivector field on  $\g^*$.
The symplectic leaves of $\g^*$ are precisely the coadjoint orbits of the corresponding Lie group $G$.
The restriction of the Poisson-Lie structure from $\g^*$ to an orbit is called the Kirillov-Kostant-Souriau (KKS) bracket.
Importance of this bracket is accounted for the fact that every $G$-homogeneous symplectic manifold is locally
isomorphic to a $G$-coadjoint orbit via the moment map.
Construction of $G$-equivariant quantization of the KKS bracket is a classic problem of the
deformation quantization theory, \cite{BFlFrLSt}.

By a quantization of  a manifold  $M$ we mean a flat $\C[[t]]$-algebra
$\A_t(M)$ with an isomorphism $\A_t(M)/t\A_t(M) \to \A(M)$, where $\A(M)$ is the function algebra on $M$.
A quantization of a  $G$-space $M$ is called equivariant if $G$ acts on $\A_t(M)$
by algebra automorphisms and that action extends the original action of $G$ on $\A(M)$.
There are various approaches to quantization of Poisson manifolds.
One of them, the $*$-product, presents a deformed multiplication in $\A(M)\tp_\C \C[[t]]$ (the tensor product
is completed in the $t$-adic topology)
as a formal $t$-series with coefficients being bi-differential operators.
Famous Fedosov's construction guarantees existence of an equivariant $*$-product on
a symplectic manifold with a $G$-invariant connection. However, Fedosov's quantization
is not given by  an explicit formula for a particular manifold.
Another approach to quantization
is to describe the  algebra $\A_t(M)$ in terms of generators and relations. While associativity
holds by the very construction and $G$-equivariance can be  easily guaranteed, the principal difficulty  is to ensure
flatness of an algebra built in such a way.


There is a universal approach to equivariant quantization of semisimple coadjoint orbits of
complex reductive Lie groups which is based on generalized Verma modules, \cite{DGS}.
By a generalized Verma module $V_{\p,\la}$ we mean the left $\U(\g)$-module $\U(\g)\tp_{\U(\p)}\C_\la$, where
$\p\subset\g$ is a parabolic subalgebra and $\C_\la$ is a one-dimensional representation of $\U(\p)$.
It is generated by a $\p$-character $\la$, which can be identified with a certain
element of $\g^*$. The approach of \cite{DGS} presents a quantized orbit as a subalgebra
of endomorphisms of a generalized Verma module.
Consider the Lie algebra $\g_t$ over $\C[t]$ which coincides with $\g[t]$ as a $\C[t]$-module and
whose Lie bracket $[\hspace{1.5pt}.\hspace{1.5pt},.\hspace{1.5pt}]_t$  extends from the bracket
of  $\g$: $[x,y]_t=t[x,y]$ for all  $x,y\in \g$.
The universal enveloping algebra
$\U(\g_t)$ is a $G$-equivariant quantization of the polynomial algebra on $\g^*$.
There is a family of isomorphisms  $\U(\g_t) \to \U(\g)$ at $t\not =0$ (this extends to a natural
embedding of $\C[t]$-algebras) such that the image of the composite map
$\U(\g_t) \to \U(\g)\to \End(\V_{\p,\la/t})$ gives an equivariant quantization of
the orbit  passing through  $\la$.

The approach of \cite{DGS} also presents quantized orbits as quotients of $\U(\g_t)$.
This fact was used in \cite{DM2} for explicit
description of quantized semisimple  $GL(n)$-orbits  in terms of generators and relations.
There was built a $\C[t]$-algebra $\A_{\mb,\mu,t}$
giving quantization of the polynomial algebra on the orbit of matrices
with eigenvalues $\mu=(\mu_1, \ldots, \mu_k)$ of multiplicities  ${\bold m}=(m_1, \ldots, m_k)$.
The algebra  $\A_{\mb,\mu,t}$ is a quotient of  $\U\bigl(\g_t\bigr)$ by
a certain ideal whose generators are given explicitly.
This ideal can be realized (not uniquely) as the annihilator of  a
generalized Verma module $V_{\p,\la/t}$, where $\p$ and $\la$ depend on $\mb$, $\mu$, and $t$.
In the present paper, we find this dependence explicitly.

Our first result is the following. We relate the two approaches to equivariant quantization
of semisimple orbits:
the quantization via generalized  Verma modules and the one in terms of generators
and relations. 

As our second result, we describe the annihilator of
a generalized Verma module $V_{\p,\la}$ for generic $\la$ as an ideal in $\U(\g)$ by presenting its set
of generators.

As an application, we classify all
finite dimensional representations of the algebras $\A_{\mb,\mu,t}$
for fixed $\mb$,  $\mu$, and $t$.
We show that they all are factored through generalized Verma modules.
For a finite dimensional representation to exist, the numbers $(\mu_i-\mu_j)/t$ should
satisfy certain integral conditions. When $\mu$ is fixed,
the dimension of representations grows as the deformation parameter $t$ tends to
zero. This process yields "fuzzyfication" for semisimple orbits of $GL(n)$.
Let us note that a reversed approach from fuzzy geometry to $*$-product quantization
was used for quantization of grassmanian spaces in \cite{DolJ}.
For applications of fuzzy geometry to theoretical physics see e.g. \cite{ARSch}.

As another application, we construct real forms  $\A_{\mb,\mu,t}$ that are compatible
with the standard real forms of $gl(n)$.

Finally, we show that special rational functions defining the central character of a quantum orbit,
 \cite{DM2}, give the Connes index for the basic quantum homogeneous  vector bundles on
 the orbit.

The paper is organized as follows.
In Sections \ref{sQGVM} and \ref{sEQ} we recall respectively the quantization
via generalized Verma modules and the quantization in terms of
generators and relations.
The correspondence  between these two approaches is established in Section \ref{sGVMQO}.
Subsection \ref{ssRQO} is devoted to finite-dimensional representations
of the algebras $\A_{\mb,\mu,t}$.
In Section \ref{sRFQO} we construct real forms on $\A_{\mb,\mu,t}$.
In Subsection \ref{ssNCCI}, we compute the non-commutative
Connes index  for the case of two-parameter quantization
of \cite{DM2}.
In Appendix we deduce an  explicit polynomial expression
for the central characters of quantum orbits.
\section{Quantization via generalized Verma modules.}
\label{sQGVM}
\begin{parag}
Consider a complex reductive Lie algebra $\g$ and let $G$ be the corresponding connected
Lie group.
The dual space $\g^*$ is a Poisson $G$-manifold
endowed with the standard Poisson-Lie structure induced by the Lie bracket
of $\g$. This bracket can be restricted to every orbit of $G$ in $\g^*$ making
it a homogeneous symplectic manifold. That symplectic structure is called
Kirillov-Kostant-Souriau bracket. An element from $\g^*$ is called
semisimple if it is the image of a semisimple element under
the $G$-equivariant isomorphism $\g\to \g^*$ via a non-degenerate
ad-invariant inner product on $\g$. An orbit is called semisimple
if it passes through a semisimple element.

Semisimple coadjoint orbits are closed affine varieties in $\g^*$. By the (polynomial) function
algebra on a subvariety  $M\subset \g^*$ we mean the algebra $\A(M)$ consisting of restrictions of
polynomial functions on $\g^*$.
In the present paper a
quantization of $M$ is a flat\footnote{By flatness of a $\C[t]$-module we mean flatness of its $t$-adic completion over $\C[[t]]$.}
$\C[t]$-algebra $\A_t(M)$ together with
an isomorphism $\kappa\colon\A_t(M)/t\A_t(M)\to \A(M)$. The quantization is called $G$-equivariant
if there is a $G$-action on $\A_t(M)$ by algebra automorphisms, an extension of the natural action
on  $\A(M)$. Such an action gives rise to a Hopf algebra action of $\U(\g)$,
$$
x\tr (ab) = (x\tr a) b+ a (x\tr b),\quad x\in \g, \quad a,b\in \A_t(M),
$$
and {\em vice versa}, so we also  call a $G$-equivariant quantization $\U(\g)$-equivariant.
\end{parag}
\begin{parag}
Let $\h$ be a Cartan subalgebra in $\g$
and $\g=\n^-\oplus\h\oplus\n^+$ a triangular decomposition
relative to $\h$. Let $\p$ be a parabolic subalgebra in $\g$
containing $\h$ and $\n^+$ and let $\l$ be its Levi factor.
The projection $\g\to \h$ along the triangular decomposition
gives rise to an embedding
$\h^*\to \g^*$.
Denote by $\c$ the center of $\l$; then $\c^*$ is
identified with the annihilator of $[\l,\l]$ in $\h^*$
via decomposition $\h=[\l,\l]\cap \h\oplus \c$.
Any $\la\in \c^*$ defines a one-dimensional representation $\C_\la$  of  $\p\supset \l$, which
extends to a representation of the universal enveloping algebra $\U(\p)$. By definition, the generalized Verma module
$\V_{\p,\la}$ is the left $\U(\g)$-module
$\U(\g)\tp_{\U(\p)} \C_\la$. When $\p$ coincides with the Borel subalgebra $\b=\h\oplus \n^+$, $\V_{\b,\la}$
is the ordinary Verma module, i.e. the maximal object in the category of $\U(\g)$-modules with the highest
weight $\la$.

Denote by $\g_t$ the Lie algebra over $\C[t]$  which coincides with $\g[t]$ as a $\C[t]$-module
and equipped with the Lie bracket $[x,y]_t=t[x,y]$ for all $x,y\in \g$.
The universal enveloping algebra $\U(\g_t)$ is a $G$-equivariant quantization
of $\g^*$.
If $\mathfrak{v}$ is a subspace in $\g$, we put $\mathfrak{v}_t =\mathfrak{v}[t]\subset\g[t]$.
Given a parabolic subalgebra $\p_t\subset \g_t$ and its character $\la$ (one dimensional representation
on $\C[t]$) we denote by $V_{\p_t,\la}$ the $\U(\g_t)$-module
$\U(\g_t)\tp_{\U(\p_t)}\C[t]$ and call it a generalized Verma module over $\U(\g_t)$.

Introduce the subset $\c^*_{reg}$
of elements in $\c^*\subset\h^*$ whose stabilizer in $\g$ is
$\l\supset \c$. We call such elements {\em regular}. The set  $\c^*_{reg}$
is a dense cone in $\c^*$. Let us fix an element
$\la\in \c^*_{reg}$.
\begin{thm}[\cite{DGS}]
\label{thm1}
The image $\A_{\p_t,\la}$ of $\U(\g_t)$ in
$\End(\V_{\p_t,\la})$  is
a $G$-equivariant quantization  of the KKS bracket
on the semisimple orbit $O_\la\subset \g^*$ passing through $\la$.
\end{thm}
Let $\J_{\p_t,\la}\subset \U(\g_t)$ denote the annihilator of the
module $\V_{\p_t,\la}$.
Theorem \ref{thm1} says that the $G$-equivariant quantization of
a semisimple orbit $O_{\la}$ can be realized as a quotient of the
algebra $\U(\g_t)$ by the ideal $\J_{\p_t,\la}$.
\end{parag}

\section{Explicit quantization for the $GL(n)$-case}
\label{sEQ}
\begin{parag}
From now on $\g$ stands for the Lie algebra $gl(n)$ of the complex general
linear group $G=GL(n)$. Let us specialize the definitions of the previous
subsection to this case. The Cartan subalgebra $\h$ is chosen to be
the subalgebra of diagonal matrices; the nilpotent subalgebras $\n^\pm$ consist of
respectively upper- and lower-triangular matrices.
Using the trace pairing, we identify the dual space $\g^*$ with $\End(\C^n)$ and $\h^*$ with the subspace of
diagonal matrices in $\End(\C^n)$; the coadjoint action of the group $G$ on $\g^*$ becomes
the similarity transformation under this identification.

Let us define the set $\{n\!:\!k\}\subset \Z^k$ of $k$-tuples
$\mb=(m_1,\ldots,m_k)$ such that $0\leq m_i$ and
$m_1+\ldots+m_k=n$; the subset $\{n\!:\!k\}_+\subset \{n\!:\!k\}$
consists of $\mb$ with all $m_i$ positive.
Let $\C^k_{reg}$ denote the subspace in $\C^k$ of $\mu=(\mu_1,\ldots,\mu_k)$ with pairwise distinct $\mu_i$.

The Levi subalgebra $\l\in \g$ consists of block diagonal matrices,
$\l=\oplus_{i=1}^k \g_i$, where $k$ is
the number of blocks and $\g_i=gl(m_i)$. We denote by $\mathrm{m}(\l)$ the $k$-tuple
$(m_1,\ldots,m_k)\in \{n\!:\!k\}_+$.
Clearly the parabolic subalgebras in $\p$ containing upper-triangular matrices are in one-to-one correspondence with elements
of $\cup_{k=1}^n\{n\!:\!k\}_+$.
Consider the decomposition $\C^n=\C^{m_1}\oplus\ldots \oplus\C^{m_k}$ and
denote by $p_i\colon\C^n\to \C^{m_i}$ the corresponding projectors.
The center of $\p$ is the linear space
$\c=\oplus_{i=1}^k \C p_i$.
Under the identification $\g^*\simeq \g$,
the dual space $\c^*$ coincides with $\c$,
and the linear isomorphism between $\C^k$ and $\c^*\subset \h^*$ is given by the
map $\mu\mapsto \sum_{i=1}^k \mu_i p_i $.
The image of the subspace $\C^k_{reg}\subset \C^k$ is the set of regular elements $\c^*_{reg}\subset \c^*$.
\end{parag}
\begin{parag}
The algebra $\U(\g_t)$ is generated by the set $\{E^i_j\}_{i,j=1}^n\subset \g_t$ of elements
satisfying the relations
\be
\label{com_rel}
\hspace{1in}E^i_j E^m_l- E^m_lE^i_j&=& t(\delta^m_j E^i_l - \delta^i_l E^m_j),\quad i,j,l,m=1,\ldots,k.
\ee
They can be arranged into a matrix, $E=||E^i_j||_{i,j=1}^n$.
The matrix $E$ can be considered as an element $E=\sum_{i,j=1}^n E^i_j\tp e^i_j\in \U(\g_t)\tp_\C \End(\C^n)$,
where $\{e^i_j\}_{i,j=1}^n$ is the standard base in $\End(\C^n)$. It is invariant with respect
to the diagonal action of $\U(\g)$ on $\U(\g_t)\tp_\C \End(\C^n)$.

We fix a multiplication in  $\End(\C^n)$ by setting it on the base by the formula
$e^i_j e^l_m = \delta_j^l e^i_m $, where $\delta_j^l$ is the Kronecker symbol.
One can consider polynomials in $E$ as an element of the algebra $\U(\g_t)\tp_\C \End(\C^n)$.
Explicitly, the $\ell$-th power $E^\ell$ is a matrix with the entries
\be
\label{monomial}
(E^\ell)^{i}_j=\sum_{\al_1,\ldots,\al_{\ell-1}=1}^n E^{\al_1}_jE^{\al_2}_{\al_1}\ldots E_{\al_{\ell-1}}^i.
\ee

We will also consider the matrix algebra $\End^\circ(\C^n)$ whose multiplication
is opposite to that of $\End(\C^n)$. The matrix $E$ may be thought of as an element from $\U(\g_t)\tp_\C \End^\circ(\C^n)$, then
its $\ell$-th power $E^{\circ \ell}$ is explicitly
\be
\label{circmonomial}
(E^{\circ \ell})^{i}_j=\sum_{\al_1,\ldots,\al_{\ell-1}=1}^n E^i_{\al_1}E_{\al_2}^{\al_1}\ldots E^{\al_{\ell-1}}_j.
\ee

In \cite{DM2}, a special series  $\Bigl\{\vartheta_\ell(\hat \mb,\mu,q^{-2},t)\Bigr\}_{\ell=0}^\infty$ of polynomials
in $\mu\in \C^k$ and $t,q^{-2}\in \C$  was introduced in connection with a two-parameter quantization on semisimple orbits.
The notation $\hat \mb$ stands for the vector $(\hat m_1,\ldots,\hat m_k)$, where the hat
denotes the  q-integers, $\hat m =\frac{1-q^{-2m}}{1-q^{-2}}$, $m\in \N$.
Here we will use the restriction $\vartheta_\ell(\mb,\mu,t)=\vartheta_\ell(\hat {\bf m},q^{-2},\mu,t)|_{q=1}$.
An explicit polynomial expression for $\vartheta_\ell({\bf m},\mu,t)$ is presented in Appendix, formula (\ref{Tr_explicit}).
\begin{definition}
\label{qOrbits}
Let $\mb=(m_1,\ldots,m_k)\in \{n\!:\!k\}_+$  and $\mu=(\mu_1,\ldots,\mu_k)\in \C^k$.
\begin{enumerate}
\item
The $\C[t]$-algebra $\A_{\mb,\mu,t}$ is a
quotient of $\U(\g_t)$ by the ideal specified by the relations
\be
\label{PR} (E-\mu_1)\ldots (E-\mu_k)&=&0,
\\
\label{TrR} \Tr\hspace{1.5pt} E^\ell&=&\vartheta_\ell(\mb,\mu,t), \quad
\ell\in 1,\ldots, k-1.
\ee
\item
The
$\C[t]$-algebra $\A^\circ_{\mb,\mu,t}$ is a
quotient of $\U(\g_t)$ by
the ideal specified by the relations
\be
\label{PR0} (E-\mu_1)\raise1.1pt\hbox{$\scriptstyle\circ$}\ldots \raise1.1pt\hbox{$\scriptstyle\circ$}(E-\mu_k)&=&0,
\\
\label{TrR0} \Tr\hspace{1.5pt} E^{\:\circ\ell} &=&\vartheta_\ell(\mb,\mu,-t), \quad
\ell\in 1,\ldots, k-1.
\ee
\end{enumerate}
\end{definition}
\noindent
The  algebras $\A_{\mb,\mu,t}$ and $\A^\circ_{\mb,\mu,t}$
are related by a certain transformation of parameters whose exact form will be presented
in Subsection \ref{A-A}, formula (\ref{link}).
We will use that relation in Subsection \ref{ssRF} concerning real forms of quantum orbits.
\begin{thm}[\cite{DM2}]
\label{final}
Given $\mu\in\C^k_{reg}$ and $\mb\in \{n\!:\!k\}_+$
the algebra  $\A_{\mb,\mu,t}$ ($\A^\circ_{\mb,\mu,t}$) is
a $G$-equivariant quantization of the semisimple orbit
of matrices with eigenvalues $\mu$ of multiplicities $\mb$.
\end{thm}
\noindent
This theorem was proven in \cite{DM2} for the algebra  $\A_{\mb,\mu,t}$. For the algebra
$\A^\circ_{\mb,\mu,t}$ the proof is analogous.
Besides, the family $\A^\circ_{\mb,\mu,t}$ can be obtain from
$\A_{\mb,-\mu,t}$ via the automorphism of the $\C$-algebra $\U(\g_t)$ extending
the map $E^i_j\to - E^i_j$, $t\to -t$.
\end{parag}
\begin{parag}
Let $S_k$ be the symmetric group of permutations of a $k$-element set.
It acts on the family of algebras $\A_{\mb,\mu,t}$ through the action on
the $k$-tuples $\mb$ and $\mu$. Denote by ${\bf l}$ the element $(1,\ldots,1)\in \C^k$.
\begin{propn}
The $\C[t]$-algebras $\A_{\mb,\mu,t}$ and $\A_{\mb',\mu',t}$ are isomorphic
if and only if
there is an element $\tau\in S_k$ and  a complex number $b$ such that
$\mb'=\tau(\mb)$ and $\mu'=\tau(\mu)+b{\bf l}$. The same holds for the family $\A^\circ_{\mb,\mu,t}$
as well.
\end{propn}
\begin{proof}
The proof is straightforward in one direction. Indeed, equations (\ref{PR})
and (\ref{TrR}) are symmetric with respect to permutations of the pairs $(m_i, \mu_i)$.
The correspondence $E^i_j\to E^i_j-b\delta^i_j$ extends to an isomorphism
$\A_{\mb,\mu,t} \to \A_{\mb,\mu+b {\bf l},t}$,
as it is seen from relations (\ref{com_rel}), (\ref{PR}), and (\ref{TrR}).

Conversely, suppose $\A_{\mb,\mu,t}$ and $\A_{\mb',\mu',t}$ are isomorphic
as $\C[t]$-algebras.
They are quantizations
of orbits characterized by eigenvalues of constituent matrices and their multiplicities.
Those orbits are
isomorphic as symplectic manifolds if and only if there is a permutation
$\tau$ and a complex number $b$ such that $\mb'=\tau(\mb)$ and $\mu'=\tau(\mu)+b{\bf l}$.
\end{proof}
\end{parag}
\begin{remark}
Given a $\C[t]$-module $V$ we will often treat it as a family of $\C$-modules.
Specialization of $\;V$ at a point $\:t=t_0$ is the $\C$-module
$V\tp_{\C[t]} \C$ corresponding to the $\C$-homomorphism $\C[t]\to \C[t]/(t-t_0)\simeq \C$.
\end{remark}
\begin{parag}[Cayley-Hamilton identity]
The elements $\Tr \:E^\ell $,
$\ell=1, \ldots, n$, generate the center $\mathcal{Z}(\g_t)$ of the algebra $\U(\g_t)$.
There is another set of generators of  $\mathcal{Z}(\g_t)$, namely, the coefficients, $\{c_i\}_{i=1}^n$,
of the "characteristic" polynomial equation
\be \label{Cayley-Hamilton}
\mathcal{P}(E)=E^n - c_1  E^{n-1}+\ldots + (-1)^n c_n=0
\ee
identically held in  $\U(\g_t)\tp_\C \End(\C^n)$.
Equation (\ref{Cayley-Hamilton}) is a non-commutative analog of the Cayley-Hamilton identity in
the classical polynomial algebra on matrices. Its existence immediately follows from
representation theory arguments and the fact that $\U(\g_t)$ is
an equivariant quantization on $\g^*$, see \cite{DM2}.

Let us define two characters of the center $\mathcal{Z}(\g_t)$ that will
appear in what follows.
\begin{enumerate}
\item
The assingment $\Tr \:E^\ell \to \vartheta_\ell(\mb,\mu,t)$, $\ell\in \N$, defines a central character, \cite{DM2},
which we denote by  $\chi_{\mb, \mu}$.
\item
Given a weight $\la\in \h^*_t$ we denote by $\chi_\la$ the central character whose
kernel lies in the annihilator of the Verma module  $V_{\b_t,\la}$.
It follows that if $\la$ is a character of a parabolic subalgebra
$\p_t\supset \b_t$, then  $\mathrm{ker}\hspace{1pt}\chi_\la$ lies in $\J_{\p_t,\la}$,
the annihilator of the generalized Verma module $V_{\p_t,\la}$.
One has $z v_0= \chi_\la(z) v_0$, where $z\in \mathcal{Z}(\g_t)$ and $v_0\in V_{\p_t,\la}$ is the
highest weight vector.
\end{enumerate}
Let $\chi$ be a character of the center $\mathcal{Z}(\g_t)$. Consider its specialization at $t\not=0$, which
is a $\C$-algebra homomorphism $\mathcal{Z}(\g_t)\to \C$.
We denote by $\Rg(\chi)$ the  set of roots of the polynomial
\be
\label{p-chi}
x^n-\chi(c_1)x^{n-1}+\ldots + (-1)^n \chi(c_n)
\ee
with complex coefficients $\chi(c_i)$.
\begin{propn}
\label{roots}
Given $\mb\in \{n\!:\!k\}_+$  and $\mu\in \C^k$ one has
\be
\Rg(\chi_{\mb, \mu})&=&\Bigl\{\mu_1,\mu_1-t,\ldots,\mu_1-(m_1-1)t;\ldots ;\mu_k,\mu_k-t,\ldots,\mu_k-(m_k-1)t\Bigr\}.
\label{rootsF}
\ee
\end{propn}
\begin{proof}
This follows from the construction of the functions $\vartheta_\ell(\mb,\mu,t)$, \cite{DM2}.
\end{proof}
\noindent
\begin{remark}
Relations (\ref{TrR}) specify an ideal in $\U(g_t)$ which is generated by the kernel of the character $\chi_{\mb, \mu}$.
Remark that for $k=n$ matrix polynomial (\ref{PR}) is obtained from (\ref{Cayley-Hamilton})
by the substitution $c_i\to \chi_{\mb, \mu}(c_i)$, as follows from (\ref{rootsF}). Thus equation (\ref{PR}) becomes superfluous, being a
consequence of (\ref{TrR}) and (\ref{Cayley-Hamilton}).
This  situation corresponds to an orbit of maximal rank, which is determined solely by a central character,  \cite{Kost}.
\end{remark}
\end{parag}
\section{Generalized Verma modules and quantum orbits}
\label{sGVMQO}
\begin{parag}
\label{sRelation}
Theorem \ref{thm1} describes quantization of  the semisimple orbit $O_\la$
as a quotient of the algebra $\U(\g_t)$ by the ideal $\J_{\p_t,\la}$.
Assuming $\la$ to be a formal function,
$\la=\la(\mu,t)$, such that $\la(\mu,0)=\mu\in \c^*_{reg}$, we obtain a family, $\A_{\p_t,\la(\mu,t)}$, of non-isomorphic
quantizations of $O_{\mu}$.
The vector $\mu$ is an element of $\C^k$ under the identification $\c^*\sim\C^k$.
On the other hand, the algebra $\A_{\mb,\mu,t}$, where $\mb=\mathrm{m}(\l)$,
is a quantization of the same orbit $O_\mu$ and it is a quotient of
$\U(\g_t)$, too.
The question is what dependence
$\la(\mu,t)$ ensures an isomorphism $\A_{\p_t,\la(\mu,t)}\simeq\A_{\mb,\mu,t}$.
\begin{lemma}[\cite{DM2}]
\label{exist}
Let $\p$ be a parabolic subalgebra with the center $\c\simeq \C^k$.
For any $\la\in \c^*_{reg}$ there exists a polynomial in one variable,
\be
\label{pp}
p(x)=x^k-\si_{1}x^{k-1}+\ldots + (-1)^k \si_{k}
\ee
with
coefficients $\si_\ell\in \C[[t]]$
such that the entries of the matrix $p(E)=||p(E)^i_j||$
lie in $\J_{\p_t,\la}$.
\end{lemma}
\noindent
For $\la\in \c^*_{reg}$,  the polynomial $p$ has only simple roots.
Rephrasing this lemma, the entries of the matrix $p(E)=||p(E)^i_j||$ annihilate the Verma module
$\V_{\p_t,\la}$.
Thus the coefficients $\si_i$ of the polynomial $p$ depend on $\la$ and $t$. At the same time they are
symmetric polynomials of its roots $\{\mu\}$.
Our goal is to find the relation between $\mu$, $\la$, and $t$.
The present section is devoted to a proof of the following
theorem.
\begin{thm}
\label{exact1}
Let $\p\subset \g$ be a parabolic subalgebra  with the Levi factor $\l$ and the center $\c\simeq \C^k$.
Put $\mb=\mathrm{m}(\l)\in \{n\!:\!k\}_+$ and take a regular element  $\la \in \c^*_{reg}\simeq \C^k_{reg}$ as
a  $\p$-character. Then
\begin{enumerate}
\item
the algebra $\A_{\mb,\mu,t} $ is isomorphic to $\A_{\p_t,\la}$ over $\C[t]$ with
\be
\label{exact_1}
\mu_i=\la_i- \sum_{\al=1}^{i-1}m_\al t,\quad i=1,\ldots,k.
\ee
\item
the algebra $\A^\circ_{\mb,\nu,t}$ is isomorphic to $\A_{\p_t,\la}$ over $\C[t]$  with
\be
\nu_i=\la_i+ \sum_{\al=i+1}^{k}m_\al t,\quad i=1,\ldots,k.
\label{exact_2}
\ee
\end{enumerate}
\end{thm}
\noindent
The rest of the section is devoted to the proof of Theorem \ref{exact1}.
We prove only statement 1; statement 2 is verified in the same manner.
\end{parag}
\begin{parag}
Consider the coefficients $\si_\ell$ of $p(x)$ in (\ref{pp}) as the elementary symmetric polynomials in $\mu\in \C^k$,
$
\si_{\ell}(\mu)=\sum_{1\leq i_1<\ldots< i_\ell\leq m}\mu_{i_1}\ldots \mu_{i_\ell}
$, $\ell=1,\ldots,k$.
First of all, to prove Theorem \ref{exact1}, we must show that the matrix entries $p(E)^i_j$
annihilate the generalized Verma module $\V_{\p_t,\la}$ provided
condition (\ref{exact_1}) holds.
Let us check the following elementary lemma.
\begin{lemma}
\label{WU=UW}
Let $W\subset\U(\g_t)$ be a submodule with respect to the adjoint representation.
The left ideal $\J_W=\U(\g_t)W$ coincides with the right ideal $W\U(\g_t)$, which is therefore a
two-sided
ideal. Then $\J_W\subset \J_{\p_t,\la}$ 
if and only if  $W$ annihilates the highest weight vector of $V_{\p_t,\la}$.
\end{lemma}
\begin{proof}
Let $\Delta$ and $\gamma$ denote  the comultiplication and antipode of the Hopf algebra
$\U(\g_t)$. In the standard Sweedler notation with implicit summation, $\Delta(x)=x_{(1)}\tp x_{(2)}$.
Let $x\in \U(\g_t)$ and $w\in W$. Then $w x= x_{(1)}\gamma\bigl(x_{(2)}\bigr)w x_{(3)}\in \U(\g_t)W$;
this proves that $\U(\g_t)W\subset W\U(\g_t)$. The reversed inclusion is checked similarly.
This  proves the first assertion of the lemma.

If $\J_W\subset \J_{\p_t,\la}$, than $W$ annihilates the highest
weight vector $v_0\in \V_{\p_t,\la}$. Conversely, suppose $W v_0=0$. An element $v\in \V_{\p_t,\la}$
can be represented as $x v_0$, where $x\in \U(\g_t)$. Hence $W v=(Wx)v_0\subset \J_W v_0=0$ and therefore
$\J_W  \subset \J_{\p_t,\la}$.
\end{proof}
\noindent
We will use Lemma \ref{WU=UW} in the situation when the module $W$ is
spanned by the entries of the polynomial $p(E)$, that is, $W=W_p=\mathrm{Span}\bigl(p(E)^i_j\bigr)$.
\end{parag}
\begin{parag}
Consider the map $\si\colon\C^n\to \C^n$, $\la_i\mapsto \si_i(\la)$,
where $\si_i$ are the elementary symmetric polynomials in $\la\in \C^n$.
The differential of this map is a matrix
$D(\la)=||D_{ij}(\la)||_{i,j=1}^n$, where $D_{ij}(\la)=\frac{\partial \si_i}{\partial \la_j}(\la)$.
\begin{lemma}
\label{wandermond}
The determinant $\mathrm{det}\: D(\la)$ is equal to
$\prod_{1\leq i<j\leq n}(\la_i-\la_j)$.
\end{lemma}
\begin{proof}
Multiply the top row of $D(\la)$, which is $||D_{1i}(\la)||=(1,\ldots,1)$, by $D(\la)_{i1}$ and subtract it
from the $i$-th row,  $i=2, \ldots, n$. This kills the entries of the first column except for
the upper one. The determinant is preserved and it equals to its  minor $\mathrm{det}||D_{ij}(\la)||_{i,j=2}^n$.
It is easy  to see, that the $j$-th column of the matrix $||D_{ij}(\la)||_{i,j=2}^n$
turns zero at $\la_{1}-\la_{j}$ and, when divided by  $\la_{1}-\la_{j}$, it
forms the $(j-1)$-th column of the $(n-1)\times(n-1)$ matrix $D(\la')$, $\la'=(\la_2,\ldots,\la_n)$.
It remains to apply induction on $n$.
\end{proof}
\end{parag}
\begin{parag} The following lemma is a key step in the proof of Theorem \ref{exact1}.
\begin{lemma}
\label{linear}
Let $\p\in \g$ be a parabolic subalgebra with the center $\c\simeq \C^k$.
There is a vector $a(\p)\in \C^k$ such that
for every $\la \in  \C^k\simeq \c^*$ the ideal  $\J_{\p_t,\la}$ contains the module $W_p$
for the polynomial $p(x)=(x-\mu_1)\ldots(x-\mu_k)$ with roots $\{\mu_i\}=\{\la_i-a_i(\p)t\}$.
\end{lemma}
\begin{proof}
First of all, let us prove this lemma assuming that $\p$ is the Borel subalgebra $\mathfrak{b}=\h+\n^+$, i.e.
for  $k=n$.
Denote by $p$ the polynomial (\ref{p-chi}) obtained from the characteristic polynomial $\mathcal{P}$
by substitution $c_i\to \chi_\la(c_i)$. Its entries annihilate the module $V_{\p_t,\la}$
because $p(E)^i_j v_0= \mathcal{P}(E)^i_j v_0=0$.
So we have $W_{p}\subset \J_{\p_t,\la}$. We will prove the statement for $\p=\mathfrak{b}$ if we show that
the roots $\{\mu\}$ of the polynomial $p$ has the form stated in the lemma.
It suffices to consider only $\la\in \C^n_{reg}$.
The coefficients $\chi_\la(c_i)$ are homogeneous polynomials
in $(\la,t)$ of degree $i$. On the other hand, they are elementary symmetric  polynomials $\si_i(\mu)$.
The map $\C^n\to \C^n$, $\mu\mapsto \si(\mu)$, is locally invertible at regular $\mu\in \C^n_{reg}$.
Hence every point in $\C^n_{reg}\times \{0\}$
has a neighborhood $U \subset \C^n_{reg}\times \C$
and an analytic function $\psi\colon U\to \C^n$ such that
$\mu=\psi(\la,t)$, $(\la,t)\in U$, are roots of the polynomial $p$.
We want to show that $\psi(\la,t)=\la - at $ for some $a=a(\b)\in \C^n$.

Observe that dilatation $E^i_j\to \frac{1}{c} E^i_j$ with $c\not=0$
extends to a $\C$-algebra isomorphism $\A_{\p_t,\la}\to \A_{\p_{ct},c\la}$.
Similarly, the dilatation $E^i_j\to \frac{1}{c} E^i_j$ extends to an isomorphism of
$\A_{\mb,\mu,t}\to \A_{\mb,c\mu,ct}$.
This implies $\psi(c\la,ct)=c\psi(\la,t)$.
Consider the $t$-expansion of the function $\psi$:
\be
\label{expand}
\psi(\la,t)=\la+\sum_{j>0} t^j\psi^{(j)}(\la),
\ee
where $\psi^{(j)}\colon \C^n\to \C^n$ are homogeneous functions  of degree $-j+1$.
Substituting $\mu= \psi(\la,t)$ to $\si_i(\mu)$
and comparing coefficients before $t^j$, we
find
$$
\sum_{l=1}^n D_{il}(\la) \psi_l^{(j)} + (\mbox{terms depending on}\> \psi^{(l)}\>\>\mbox{with}\>\> l<j)\in \C^n[\la].
$$
Using induction on $j$ and Lemma \ref{wandermond}, we see that $\psi^{(j)}_i$ are rational functions in $\la$
maybe having poles only at $\la_i=\la_l$, $i\not =l$.
But $\psi^{(j)}$ are bounded at  $\la_i=\la_l$ since $\mu=\psi(\la,t)$  are roots of the polynomial
$p$ with the coefficients $\chi_\la(c_i)$ being regular functions in $(\la,t)$.
Therefore $\psi^{(j)}\in \C^n[\la]$ and, taking into account their homogeneity degree $-j+1$,
all they are zero except for the  term $\psi^{(1)}=-a(\b)$.
Thus we have proven the statement for $\p=\b$.
As a corollary, we obtained $\Rg(\chi_\la)=\{\la-a(\b)t\}$.

Let us consider the situation of a general parabolic subalgebra $\p\supset \b$.
Take $\la\in \c^*_{reg}$ and denote $\tilde\la$ its image under the canonical embedding in $\c^*\subset\h^*$.
By Lemma \ref{exist}, there is a polynomial $p$ of degree $k$ whose entries lie in $\J_{\p_t,\la}$.
Let $\mu$ be its roots. The central character $\chi_{\mb,\mu}$ coincides
with $\chi_{\tilde \la}$, therefore $\{\mu\}\subset \{\tilde \la-a(\b)t\}$
by Proposition \ref{roots}. This proves the lemma
when $\la\in \c^*_{reg}$. But the condition $W_p v_0=0$
determines the coefficients of $p$ as rational functions of $(\la,t)$.
We have already proven that they are in fact polynomials, being the elementary symmetric
polynomials in the roots. Therefore
the polynomial $p$ with the property $W_p \subset \J_{\p_t,\la}$ does exist for all $\la$.
Its roots $\{\mu\}$ are contained in  $\{\tilde \la-a(\b)t\}$, so
$\mu$ is related to $\la$ as stated in the lemma.
\end{proof}
\begin{remark}
\label{unique}
The module $W_p$ is unique for regular $\la$ and sufficiently small $t$. Indeed,
the condition $W_p v_0=0$ gives rise to a system of equations on the polynomial $p$. It goes over to
the system
$p(\la_i)=0$, $i=1\ldots k$,  at $t=0$. This uniquely determines $p$ up to a factor at
$\la\in \c^*_{reg}$ and $t=0$, hence $p$ is unique for $\la\in \c^*_{reg}$ and small $t\not=0$.
\end{remark}
\noindent
To finish the proof of Theorem \ref{exact1}, it remains to determine the vector $a(\p)\in \C^k$.
\end{parag}
\begin{parag}[Proof of Theorem \ref{exact1}]
\label{ssAGVM}
Let $\p$ be the parabolic subalgebra with the center $\c\simeq \C^k$ and put $\mb=m(\l)\in \{n\!:\!k\}$.
In this subsection we show that the $j$-th coordinate of the vector $a(\p)\in \C^k$ from Lemma \ref{linear}
is equal to $\sum_{i=1}^{j-1}m_i$. This will prove Theorem \ref{exact1}.
Let us first consider the case $k=2$.
Suppose $\la\in \c^*_{reg}$.
Only two equations are independent in
the system $p(E)^i_i v_0=0$, $i=1, \ldots,n$,
 say for $i=1$ and $i={m_1+1}$:
\begin{equation}
\begin{array}{lcc}
\la_1^2-(\mu_1+\mu_2)\la_1+\mu_1\mu_2&=&0,\nn\\
\la_2^2-(\mu_1+\mu_2+tm_1)\la_2+\mu_1\mu_2+t m_1\la_1&=&0.
\end{array}
\nn
\end{equation}
This system has a solution $(\la_1,\la_2)=(\mu_1,\mu_2+tm_1)$  satisfying the condition
$\la=\mu$ at $t=0$. Thus we have proven Theorem \ref{exact1} for $k=2$.

We deduce  the case $k=n$ from the studied case $k=2$.
Recall that given an element $\tilde \la\in \h^*$ the set $\mathfrak{R}(\chi_{\tilde \la})$ is
equal to
$\bigl\{\tilde\la_i-a_i(\b)t\bigr\}^n_{i=1}$.
Putting $\tilde\la_i=\la_1$ for $i=1,\ldots, m_1$ and $\tilde\la_i=\la_2$ for $i=m_1+1,\ldots ,n$,
we obtain a central character $\chi_{\mb,\mu}$ associated with the algebra $\A_{\mb,\mu,t}$,
for
$\mb\in \{n\!:\!2\}_+$, and $\mu\in \C^2$.
We have already shown that $\mu_1=\la_1$ and $\mu_2=\la_2-tm_1$.
Thus we have, by Proposition \ref{roots},
$$
\mathfrak{R}(\chi_{\mb,\mu})=\bigl\{\la_1-(i-1)t\}_{i=1}^{m_1}\cup \{\la_2-(m_1+i-1)t\bigr\}_{i=1}^{m_2}.
$$
On the other hand,
$$
\Rg(\chi_{\tilde \la})=
\bigl\{\la_1-a_i(\b)t\bigr\}_{i=1}^{m_1}\cup \bigl\{\la_2-a_{m_1+i}(\b)t\bigr\}_{i=1}^{m_2}.
$$
The sets $\Rg(\chi_{\mb,\mu})$ and $\Rg(\chi_{\tilde \la})$ coincide since $\chi_{\mb,\mu}= \chi_{\tilde \la}$.
Let us put $t=1$ and take $\la_1$ and $\la_2$ to be positive and
negative real numbers, respectively. We can assume them big enough in their absolute values so as to
ensure coincidence of the subsets
$\bigl\{\la_1-(i-1)\bigr\}_{i=1}^{m_1}\subset \mathfrak{R}(\chi_{\mb,\mu})$ and
$\bigl\{\la_1-a_i(\b)\bigr\}_{i=1}^{m_1}\subset\Rg(\chi_{\tilde \la})$.
Using this argument, we can apply induction on $m_1$ and prove the case $k=n$ of Theorem \ref{exact1}.
This gives $a_i(\b)=(i-1)$, $i=1,\ldots,n$.

It remains to consider the case $2<k<n$. Let us set $t=1$ and assume that the coordinates
of $\la\in \c^*_{reg}\simeq \C^k$ are real and form a strictly decreasing sequence, $\la_j>\la_{j+1}$.
The vector $\la \in \C^k$ corresponds to the vector
$\tilde\la\in \h^*\simeq \C^n$ via the canonical embedding $\C^k\simeq\c^*\to \h^*\simeq \C^n$;
the coordinates of $\tilde \la$
are $(\la_1,\ldots,\la_1;\ldots;\la_k,\ldots,\la_k)$,
where each $\la_i$ is taken $m_i$ times.
The set $\mathfrak{R}(\chi_{\tilde \la})$ is a union $\cup_{j=1}^k \mathfrak{R}_j$ of non-intersecting subsets
\be
\label{R1}
\mathfrak{R}_j=\bigl\{\la_j-(\sum_{l=1}^{j-1}m_l+ i-1)\bigr\}_{i=1}^{m_j}.
\ee
We have $\min \mathfrak{R}_j > \max \mathfrak{R}_{j+1}$.
On the other hand, we know from Proposition \ref{roots} and Lemma \ref{linear} that $\mathfrak{R}(\chi_{\tilde \la})$
is a union the of subsets
\be
\label{R2}
\mathfrak{R}'_j=\bigl\{\mu_j-(i-1)\}_{i=1}^{m_j}=\{\la_j-(a_j(\p)+i-1)\bigr\}_{i=1}^{m_j}.
\ee
We can chose such $\la_j$  that
$\mathfrak{R}'_j$ do not intersect and, moreover,
$\min \mathfrak{R}'_j > \max \mathfrak{R}'_{j+1}$
(for that to be true, it is enough to assume $\la_j-\la_{j+1}>a_j(\p)-a_{j+1}(\p)+m_j-1$).
Taking into account $\#\hspace{3pt}\mathfrak{R}'_j=m_j=\#\hspace{1pt}\mathfrak{R}_j$ we conclude that
$\mathfrak{R}'_j=\mathfrak{R}_j$ for all $j$. Comparing (\ref{R1}) and (\ref{R2}), we find $a(\p)$  and
thus prove Theorem \ref{exact1}.
\end{parag}
\begin{remark}
\label{Verma}
As follows from Lemma \ref{linear}, the identity map $\g_t\to \g_t$ extends to an epimorphism
$\A_{\mb,\mu,t}\to \U(\g_t)/\J_{\p_t,\la}$ of $\C$-algebras if the pair $(\p_t,\la)$ is
related to the pair $(\mb,\mu)$ as in (\ref{exact_1}). For regular $\la$ this
map is an isomorphism of $\C[t]$-algebras.
It is known,  \cite{J}, that for almost all $\la\in \c^*$ the
annihilator of  $V_{\p,\la}$ is generated by a copy of adjoint representation
and the kernel of a central character. Theorem \ref{exact1} gives an explicit description
of the annihilator.
\end{remark}
\begin{parag}[Representations of quantum orbits]
\label{ssRQO}
This subsection is devoted to finite dimensional representations of the $\C$-algebras $\A_{\mb,\mu,t}$.
Our principal tool is the established correspondence between the family  $\A_{\mb,\mu,t}$ and generalized
Verma modules. Note that this correspondence is not one-to-one.
The algebra $\A_{\mb,\mu,t}$ can be realized as the algebra $\A_{\p_t,\la}$ for some $\p_t$ and $\la$
in different ways.
This freedom comes out from permutations $(\mb,\mu)\to \bigl(\tau(\mb),\tau(\mu)\bigr)$,
$\tau\in S_k$, leaving $\A_{\mb,\mu,t}$ invariant but changing the generalized Verma modules.
The situation is exactly the same as in the case $\p=\b$ when the Weyl group transitively
acts on the set of the ordinary Verma modules
with isomorphic annihilators.

Any representation of $\A_{\mb,\mu,t}$ is at the same time a representation of $\U(\g_t)$.
Therefore to describe representations of $\A_{\mb,\mu,t}$ it is necessary and sufficient  to
determine those representation of $\U(\g_t)$ that are factored through the ideal specifying
$\A_{\mb,\mu,t}$. We do it for finite dimensional representations.
\begin{propn}
The quantum orbit $\A_{\mb,\mu,t}$   has a finite dimensional representation
if and only if there is
an element $\tau\in  S_k$ such that $(\mu'_{i}-\mu'_{i+1})/t-m'_{i}$, where
$\mu'=\tau(\mu)$ and $m'=\tau(m)$, are non-negative integers for $i=1,\ldots,k-1$.
If exists, such a representation is unique and it is factored through a generalized Verma module.
\end{propn}
\begin{proof}
As was already mentioned in the proof of Lemma \ref{linear}, the correspondence
$E^i_j \to \frac{1}{c}E^i_j$, where $c\not =0$, extends to an isomorphism between the  algebras
$\A_{\mb,\mu,t}$ and $\A_{\mb,c\mu,ct}$. Therefore we can assume $t=1$ in the proof.

Suppose there is such element $\tau \in S_k$ as stated in the proposition.
Take the parabolic subalgebra $\p$ with the Levi factor $\l$
such that $\mathrm{m}(\l)=\mb'$. Consider the $\U(\g)$-module $V_{\p, \la}$, where $\la$ is related
to $\mu'$ by formula (\ref{exact_1}), where one should set $t=1$.
By assumption, the numbers $\la_i=\mu'_{i}+ \sum _{\al=1}^{i-1}m'_{\al}$
define a dominant  integral weight of $sl(n)\subset \g$, so the
module $V_{\p, {\bf \la}}$ is projected onto a finite dimensional
$\U(\g)$-module $W$ with the highest weight $\la$.
The homomorphism of modules induces a representation of  $\A_{\mb,\mu,1}\simeq \A_{\mb',\mu',1}$
on $W$.

Conversely, suppose  $W$ is a finite dimensional module over  $\A_{\mb,\mu,1}$.
Then it is a module over  $\U(\g)$ as well.
There is a cyclic $\U(\g)$-submodule in $W$ with the highest weight $\la$.
We can  think that this submodule coincides with $W$.
The highest weight gives rise to a central character, $\chi_\la$.
On the other hand, one has $\chi_\la=\chi_{\mb,\mu}$,
and the set $\mathfrak{R}(\chi_\la)=\bigl\{\la_j-(j-1)\bigr\}_{j=1}^{n}$
contains  $\{\mu_i\}_{i=1}^k$ as a subset, by Proposition \ref{roots}.
Further, $W$ is a finite dimensional $\U(\g)$-module with the highest
weight $\la$, therefore $\la_j-\la_{j+1}$ are non-negative integers.
The elements of $\mathfrak{R}(\chi_\la)$ are ordered  by their real components; then
the numbers $\la_j-(j-1)t\in\mathfrak{R}(\chi_\la)$, $j=1,\ldots,n$, form a strictly decreasing sequence.
The subset $\{\mu\}\subset \mathfrak{R}(\chi_\la)$ is ordered by
inclusion, hence there is a permutation $\tau$ such that
$\mu'_i>\mu'_{i+1}$, $\mu'=\tau(\mu)$. The permutation $\tau$ satisfies the
condition of the theorem. Uniqueness of the module $W$ follows from
the ordering on $\{\mu'\}$.

Let us show that the module $W$ is a quotient of a generalized
Verma module associated with $\A_{\mb,\mu,1}$.
The elements $\la_j-(j-1)$ belonging to the interval
$\bigl[\mu'_i,\mu'_{i+1}\bigr)\subset \mathfrak{R}(\chi_\la)$
form an arithmetic progression with the initial term $\mu'_i$ and decrement $1$, as follows from
(\ref{rootsF}). Thus
$\la_j$ is stable within this interval. Therefore  $\la$ defines a character of the parabolic subalgebra
$\p$ with Levi factor $\l$ such that $\mathrm{m}(\l)=\mb'$.
Let $\rho$ denote the homomorphism $\U(\g)\to \End(W)$.
As $\dim W <\infty$,  for all $x\in \l$ one has $\rho(x)w_0=\la(x)w_0$, where
$w_0\in W$ is the highest weight vector.
Therefore the subspace $\C w_0\subset W$ forms a one dimensional $\U(\p)$-module $\C_\la$.
The  map  $\U(\g)\tp_\C \C \to W$, $x \tp 1 \mapsto \rho(x)w_0$, is a homomorphism of left $\U(\g)$-modules, which is obviously factored through the
generalized Verma module $\U(\g)\tp_{\U(\p) }\C_{\la}$.
\end{proof}
\end{parag}
\section{Real forms of quantum orbits}
As an application of the established relation between the algebras $\A_{\mb,\mu,t}$ ($\A^\circ_{\mb,\mu,t}$) and
the generalized Verma modules we construct real forms of quantum orbits.
\label{sRFQO}
\begin{parag}
An anti-linear endomorphism of a complex vector space $V$ is
an additive map $f\colon V\to V$ satisfying $f(\bt a)=\bar \bt f(a)$
for any $a\in V$ and $\bt\in\C$; the bar stands for the complex conjugation.
When $V$ is a $\C[t]$-module, we assume $f(t)=t$.
\begin{definition}
A $*$-structure on an associative algebra $\A$ over $\C$ or $\C[t]$  is
an anti-linear involution $*\colon \A\to \A$ that is an anti-automorphism
with respect to the  multiplication,
$(ab)^*=b^*a^*$.
\end{definition}
\noindent
Let $\Ha$ be a Hopf algebra over $\C$ or $\C[t]$ with the multiplication $m$,  comultiplication
$\Delta$, and antipode $\gamma$, \cite{ChPr}.
\begin{definition}
A real form on $\Ha$   is an anti-linear involution such that
\be
\theta\circ m =m\circ(\theta\tp\theta), \quad \tau\circ\Delta\circ \theta =(\theta\tp\theta)\circ\Delta,
\ee
where $\tau$ is the flip on $\Ha\tp \Ha$.
\end{definition}
\noindent
When $\Ha$ is the universal enveloping algebra $\U(\g)$ of a complex Lie algebra $\g$,
the involution restricts to $\g$. The set of  $\theta$-fixed points in  $\U(\g)$ is
the universal enveloping algebra $\U(\g_\R)$ over $\R$ of the Lie algebra
$\g_\R$, a real form of $\g$.
\begin{propn}
Let $\theta$ be a real form of a Hopf algebra $\Ha$ with invertible antipode $\gm$.
Then
$\gm\circ \theta \circ \gm = \theta$.
\end{propn}
\begin{proof}
Consider the map $\theta \circ \gm\circ\theta \colon \Ha\to \Ha$. It satisfies
the axioms of antipode for the opposite comultiplication and therefore
equals to $\gm^{-1}$, due to uniqueness of the antipode.
\end{proof}
\noindent
This proposition implies that the composition of $\theta$ with any odd power of $\gm$ makes $\Ha$ an involitive algebra.
We are going to define involutions on $\Ha$-module algebras that will be compatible with real forms on $\Ha$
in the sense of the following definition.
\begin{definition}
Let $\theta$ be a real form of a Hopf algebra $\Ha$.
A real form of an $\Ha$-module algebra $\A$ is a $*$-structure on $\A$ such that
\be
(x\tr a)^*&=& \theta(x)\tr a^*,\quad x\in \Ha, \quad a\in \A.
\ee
\end{definition}
\begin{example}
\label{*-theta}
Let $\Ha$ be a Hopf algebra equipped with  a real form $\theta$.
Consider $\Ha$ as a left module algebra over itself with respect
to the adjoint action $x\tr y= x_{(1)}y \gm\left(x_{(2)}\right)$.
Define a $*$-structure on $\Ha$ by the involution $\gm \circ \theta$.
Then, $*$ is a real form of the self-adjoint module algebra $\Ha$.
Indeed, one has
$$
(x\tr y)^* = (x_{(1)}y \gm(x_{(2)}))^*=(\gm\circ\theta\circ\gm)(x_{(2)}) y^* (\gm\circ\theta)(x_{(1)})
=\theta(x_{(2)}) y^* (\gm\circ\theta)(x_{(1)}),
$$
which is equal to $\theta(x)\tr y^*$.
\end{example}
\end{parag}
\begin{parag}
\label{A-A}
Before proceeding with real forms as applied to $G$-equivariant quantization,
we study a relation between the two algebras $\A_{\mb,\mu,t}$ and  $\A^\circ_{\mb,\mu,t}$
from Definition \ref{qOrbits}.
\begin{propn}
The identity map $\g_t\to \g_t$ induces an isomorphism of  algebras  $\A_{\mb,\mu,t}$ and  $\A^\circ_{\mb,\nu,t}$
over $\C$,
provided
\be
\label{link}
\nu_i=\mu_i+(n-m_i)t, \quad i=1,\ldots, k.
\ee
\end{propn}
\begin{proof}
Both algebras are quotients of $\U(\g_t)$ by the ideals generated by certain copies of
the adjoint and trivial submodules. Denote them by $W_\mu$ and $W^\circ_\nu$ correspondingly.
It is sufficient to show that $W_\mu$ and $W^\circ_\nu$ are related by the transformation
of parameters (\ref{link}).
With $\mb$ given, both modules are determined by their adjoint components,
generated by the matrix coefficients of polynomials in $E$.
Every monomial (\ref{monomial}) is expressed as a linear combination of
monomials (\ref{circmonomial}) with coefficients being polynomial
functions in $t$. A polynomial (\ref{PR})
can be represented as a polynomial in the sense of (\ref{circmonomial}),
and the  coefficients of the latter will be polynomials in $\mu$ and $t$.
Therefore it suffices  to prove the proposition
only for generic $\mu$.

Take the parabolic subalgebra $\p\supset \l$ such that
$m(\l)=\mb$. Assuming $\mu$ regular  and $t$ small, there is
a generalized Verma module such that
$W_\mu\subset \J_{\p_t,\la}$. Further, there is a unique, as emphasized in Remark \ref{unique},
module $W^\circ_\nu$ such that $\J_{\p_t,\la}\supset W^\circ_\nu$.
Hence $W_\mu=W^\circ_\nu$, and the
parameters $\mu$, $\la$, and $\nu$ are related by transformations (\ref{exact_1}) and  (\ref{exact_2}).
This proves the statement for regular $\mu$ and small $t$ and therefore for all
$\mu$ and $t$.
\end{proof}
\end{parag}
\begin{parag}
It is known, \cite{Kn}, that all real forms of the Lie algebra $gl(n,\C)$ are
isomorphic to
$gl(n,\R)$, $u^\bigstar(2m)$ when $n=2m$ is even, and $u(r,s)$ with $r+s=n$, $r\geq s\geq 0$.
When $s=0$, $u(r,s)$ turns into the compact real form $u(n)$.
The corresponding involutions are defined on generators $\{E^i_j\}_{i,j=1}^k$ as
$$
\theta_{gl(n,\R)}(E)=E, \quad  \theta_{u(r,s)}(E)=-J_{(r,s)}E' J^{-1}_{(r,s)},
\quad  \theta_{u^\bigstar(2m)}(E)=J_{m} E J^{-1}_{m},
$$
where  $J_{(r,s)}=\sum_{i=1}^r e^i_i -\sum_{i=r+1}^n e^i_i $,
$J_{m}=\sum_{i=1}^{m} (e^i_{i+m} -e^{i+m}_i)$,
and the prime stands for the matrix transposition.
\end{parag}
\begin{parag}
\label{ssRF}
Take $\mu\in \C^k$, $\mb\in \{n\!:\!k\}$ and denote by $S^\mb_k$ the stabilizer of $\mb$
in the symmetric group $S_k$.
\begin{propn}
\label{su}
The map $E\to J_{(r,s)}E' J^{-1}_{(r,s)}$ extends to
a $\theta_{u(r,s)}$-compatible $*$-real form of the algebra $\A_{\mb,\mu,t}$,
provided
$\bar \mu = \tau(\mu)$ for some $\tau\in S^\mb_k$.
\end{propn}
\begin{proof}
Consider $\U(\g)$ as a module algebra over itself with respect to
the adjoint representation. Define a $\theta_{u(r,s)}$-compatible
$*$-real form on  $\U(\g)$ as in Example \ref{*-theta}.
It is naturally extended to a  $\theta_{u(r,s)}$-compatible
$*$-real form on $\U(\g)[t]$ considered as a  $\U(\g)$-module
algebra. This real form is restricted to the $\U(\g)$-equivariant
embedding $\U(\g_t)\to \U(\g)[t]$.
On the generators $\{E^i_j\}\subset \g_t$, it is defined by the map $E\to J_{(r,s)}E' J^{-1}_{(r,s)}$.
Let us prove that it induces an isomorphism $(\A_{\mb,\mu,t})^* \simeq \A_{\mb,\bar\mu, t}$.
Setting  $J=J_{r,s}$  we have $J^2=1$ and $J'=J$; hence for $m>0$
\be
(E^\ell)^i_j&=&
\sum_{\al_1,\ldots, \al_{\ell-1}}
E_j^{\al_1}E_{\al_1}^{\al_2}\ldots E_{\al_{\ell-1}}^{i}
\stackrel{*}{\to }\sum_{\al_0,\ldots, \al_\ell}
J_i^{\al_\ell}E^{\al_{\ell-1}}_{\al_\ell}\ldots E^{\al_1}_{\al_2} E^i_{\al_1} J^j_{\al_0}.
\nn
\ee
In the concise matrix form this reads $(E^\ell)^*= (J E^\ell J)'$.
If $p(x)$ is a polynomial in one variable and $\bar p(x)$ is obtained from
$p(x)$ by the complex conjugation of  its coefficients, then
$p(E)\stackrel{*}{\to}\bigl(J\bar p(E)J\bigr)'$.
Relations (\ref{TrR}) go over to
$
\Tr\hspace{1.5pt} E^\ell =\vartheta_\ell(\mb,\bar \mu, t)$,
$\ell\in \N$.
Therefore $(\A_{\mb,\mu,t})^* \simeq \A_{\mb,\bar\mu, t}=\A_{\mb,\mu,t}$ if
 $\bar \mu = \tau(\mu)$ for some $\tau\in S^\mb_k$.
\end{proof}
\begin{propn}
Suppose $-\bar \mu =\tau(\mu-t\mb) +t  n {\bf l}\in \C^k$ for some $\tau\in S^\mb_k$.
Then
\begin{enumerate}
\item
the map $E\to E$ extends to
a $\theta_{gl(n,\R)}$-compatible $*$-real form of $\A_{\mb,\mu,t}$
\item
the map $E\to -J_{m} E J^{-1}_{m}$ extends to
a $\theta_{u^\bigstar(2m)}$-compatible $*$-real form of $\A_{\mb,\mu,t}$
\end{enumerate}
\end{propn}
\begin{proof}
Let $J$ denote either the unit matrix or $J_{m}=\sum_{i=1}^{m} (e^i_{i+m} -e^{i+m}_i)$.
We will consider the two cases simultaneously.
Define a $*$-real form on the algebra $\U(\g_t)$ in the same way as
in the proof of Proposition \ref{su}. It is compatible with
the corresponding real form of the Hopf algebra $\U(\g)$.
On the generators $\{E^i_j\}\subset \g_t$, it is defined by the map $E\to -J E J^{-1}$.
We have
\be
(E^\ell)^i_j&=&
\sum_{\al_1,\ldots, \al_{\ell-1}}
E_j^{\al_1}E_{\al_1}^{\al_2}\ldots E_{\al_{\ell-1}}^{i}
\stackrel{*}{\to }\sum_{\al_0,\ldots, \al_\ell}
(-1)^\ell J^i_{\al_\ell}E_{\al_{\ell-1}}^{\al_\ell}\ldots E_{\al_1}^{\al_2} E_{\al_0}^{\al_1}\bigl(J^{-1}\bigr)_j^{\al_0}.
\nn
\ee
Thus every matrix monomial $E^\ell$ is transformed into the matrix $(-1)^\ell J E^{\circ \ell} J^{-1}$; therefore
$(\A_{\mb,\mu,t})^*\simeq\A^\circ_{\mb,-\bar \mu,t}$, with respect
to either $sl(n,\R)$- or $u^\bigstar(2m)$-involutions.
Now the rest of the proof follows from (\ref{link}).
\end{proof}
\end{parag}
\section{Non-commutative Connes index}
\label{sNCCI}
\begin{parag}
\label{Polynomials}
The purpose of the section is to give an interpretation of certain rational functions
arising in our theory of quantum orbits, \cite{DM2}. We show that those functions give
the non-commutative  Connes index, \cite{C}, of basic projecitve modules over quantum orbits.
We will consider the two-parameter $\U_q(\g)$-equivariant quantization of \cite{DM2} (see also \cite{DM1})
including the $G$-equivariant quantization as the limit case $q\to1 $.
A two-parameter quantum orbit is a quotient of the so called modified reflection equation
algebra, \cite{KSkl,IP,DM3}, $\L_{q,t}$, which itself is a  $\U_q(\g)$-equivariant quantization
of the polynomial algebra on $\g^*$.
The algebra  $\L_{q,t}$ is generated by the elements $\{L^i_j\}_{i,j=1}^n$ subject to
certain quadratic-linear relations turning to
commutation relations (\ref{com_rel}), where $E^i_j=\lim_{q\to1}{L^i_j}$.
There is a $\U_q(\g)$-equivariant generalization, $\Tr_q$, of the trace operation such that
$\Tr_q \:L^\ell$, $\ell\in \N$, belong to the center of  $\L_{q,t}$.
For details the reader is referred to \cite{DM2} and references therein.

The functions $\vartheta_\ell(\mb, \mu, t)$ in the right-hand side of (\ref{TrR}) are the specialization of
certain functions $\vartheta_\ell(\hat \mb, \mu, q^{-2},t)$
introduced in \cite{DM2} (see also Appendix).
Here,  $\hat \mb = (\hat m_1,\ldots, \hat m_k)$ and $\hat m=\frac{1-q^{-2m}}{1-q^{-2}}$ for $m\in \N$.
Also,  we introduced in \cite{DM2} rational functions $C_j(\hat \mb,\mu,q^{-2},t)$ satisfying
relation (\ref{qtConnesI}).
\begin{thm}[\cite{DM2}]
\label{final_q} For any  $\mb\in
\{n\!:\!k\}$ and $\mu\in \C^k_{reg}$, the quotient $\A_{\hat\mb,\mu,q,t}$ of the algebra $\L_{q,t}$ by
the relations
\be \label{PRq} (L-\mu_1)\ldots (L-\mu_k)&=&0,
\\
\label{TrRq}
\Tr_q(L^\ell)&=&\vartheta_\ell(\hat\mb,q^{-2},\mu,t), \quad
\ell=1,\ldots,k-1,
\ee
is a $\U_q(\g)$-equivariant quantization on the orbit
of semisimple matrices with eigenvalues $\mu$ of
multiplicities $\mb$.
\end{thm}
\end{parag}
\begin{parag}
\label{ssNCCI}
In the classical geometry, a vector bundle over
a manifold $M$ may be given by an idempotent $\pi$ of the algebra $\A(M)\tp_\C \End(V)$,
where $V$ is a finite dimensional vector space. It forms a projective $\A(M)$-module
$\pi \bigl(\A(M)\tp_\C V\bigr)$.
If $M$ is a
$G$-manifold and $V$ a $G$-module, a $G$-equivariant bundle corresponds to an invariant
idempotent.
For every semisimple coadjoint orbit $O_{\mu}$ in $\End(\C^n)$ of rank $k-1$ there are $k$
projector-valued functions $\pi_i\colon O_{\mu} \to \End(\C^n)$.
At a point $A\in O_{\mu}$, they commute with $A$ and map the linear space
$\C^n$ onto the  $A$-eigenspaces.
These basic vector bundles generate the Grothendieck ring
of equivariant vector bundles over $O_{\mu}$.

In the two-parameter quantization setting an idempotent $\pi$
from $\A_{q,t}(M)\tp_\C \End(V)$ defines the projective
$\A_{q,t}$-module $\pi \bigl(\A_{q,t}(M)\tp_\C V\bigr)$, which we
consider as a quantized vector bundle over $M$.
Let us  built the idempotents for the quantized basic vector bundles over the orbits.
\begin{propn}
Let $L=||L^i_j||_{i,j=1}^k$ be the matrix of generators of the algebra $\A_{\hat\mb,\mu,q,t}$ and  suppose
$\mu\in \C^k_{reg}$.
The elements
\be
\label{projectors}
\pi_j(L)&=&\prod_{i=1,\ldots, k\atop i\not= j}\frac{L-\mu_i}{\mu_j-\mu_i} \in \A_{\hat\mb,\mu,q,t}\tp_\C \End(\C^n), \quad j=1,\ldots,k,
\ee
are invariant projectors.
One has $\Tr_q\hspace{1.5pt} \pi_j(L)= C_j(\hat \mb,\mu,q^{-2},t)$, $j=1,\ldots,k$, see (\ref{qtConnesI}).
\end{propn}
\begin{proof}
The matrix $L\in \A_{\hat\mb,\mu,q,t}\tp_\C \End(\C^n)$ satisfies (\ref{PRq}).
This implies that the elements (\ref{projectors}) are orthogonal idempotents,
and $L=\sum_{j=1}^k \mu_j \:\pi_j(L)$. Therefore, every positive integer power of the
matrix $L$ has the decomposition over the basis $\pi_j(L)$, $j=1,\ldots,k$:
\be
L^\ell=\sum_{j=1}^k \mu^\ell_j \:\pi_j(L), \quad\ell=\N.
\ee
Taking trace on
the both sides of this equation and using condition (\ref{TrRq}) extended
for all $\ell\in \N$ (see \cite{DM2})
and representation (\ref{qtConnesI}), we get
$\sum_{j=1}^k \mu^\ell_j \:\Tr_q \hspace{1.5pt}\pi_j(L)=\sum_{j=1}^k \mu^\ell_j\: C_j(\hat \mb,\mu,q^{-2},t)$
for all $\ell=\N$. Since the numbers $\{\mu_i\}$ are pairwise distinct, this proves the statement.
\end{proof}
The quantities $\Tr_q \hspace{1.5pt}\pi_j(L)$ are central elements
of  $\A_{\hat\mb,\mu,q,t}$. Therefore we can take the function
$C_j(\hat \mb,\mu,q^{-2},t)$ as the "universal" Connes index of the family of quantum vector bundles
$\pi_j \bigl(\A_{\hat\mb,\mu,q,t}\tp_\C \C^n\bigr)$.

Note that the problem of $\U_q(\g)$-equivariant quantization of vector bundles on orbits
as one- and two-sided projective modules over quantized function algebras was considered in \cite{D}.
An interesting problem is to construct the quantum bundles explicitly, in terms of projectors.
This problem is solved for  non-commutative sphere $\mathbb{S}^2_q$ in \cite{GLS}.
\end{parag}

\vspace{0.8cm}
\noindent
{\bf \Large Appendix}

\vspace{10pt}
\noindent
In this subsection, we study the functions $\vartheta_\ell$ entering
the right-hand side of (\ref{TrR}). They were introduced in \cite{DM2}
as specialization at $q=1$ of certain functions participating
in the two-parameter quantization of orbits.
Assuming $\nu,\mu\in \C^k$ and $\omega\in \C$ consider the functions
\be
\label{tr_explicit}
\tilde\vartheta_\ell(\nu,\mu,\omega)&=&\sum_{s=1}^{k} \omega^{s-1}
\sum_{1\leq j_1<\ldots <j_s \leq k}\nu_{j_1}\ldots \nu_{j_s}
\sum_{{\bf d}\in \{\ell:s\}}\mu^{d_1}_{j_1}\ldots \mu^{d_s}_{j_s}, \quad \ell=1,2,\ldots
\ee
 They satisfy the recurrent relation
\be
\label{n_n-1}
\tilde\vartheta_\ell(\nu,\mu,\omega) =
\tilde\vartheta_\ell(\nu',\mu',\omega)
+\nu_k \omega\sum_{i=1}^{m-1}
\tilde\vartheta_{\ell-i}(\nu',\mu',\omega)\mu_k^i + \nu_k \mu_k^m,
\ee
where
$\mu'=(\mu_1,\ldots,\mu_{k-1})$ and $\nu'=(\nu_1,\ldots,\nu_{k-1})$.
Using this relation, it is easy to prove by induction on $k$ that
\be
\label{thetas}
\tilde\vartheta_\ell(\nu,\mu,\omega)=\sum_{j=1}^k \mu_j^\ell\:\tilde C_j(\nu,\mu,\omega),
\ee
where
\be
\label{ConnesI}
\tilde C_j(\nu,\mu,\omega)&=& \nu_j+\nu_j\sum_{\ell=1}^{k-1}
\omega^{\ell} \sum_{1\leq i_1<\ldots <i_\ell\leq k\atop i_1,\ldots, i_\ell\not=j} \frac{\nu_{i_1}
\mu_{i_1}}{\mu_{j}-\mu_{i_1}}\ldots \frac{\nu_{i_\ell} \mu_{i_\ell}}{\mu_{j}-\mu_{i_\ell}} \:,
\ee
If one puts, by definition, $\tilde\vartheta_0(\nu,\mu,\omega) =\nu_\omega$,
where $\nu_\omega$ is determined by the equation $(1-\omega\nu_\omega)=\prod_{i=1}^k(1-\omega\nu_i)$,
then representation (\ref{thetas}) is valid for $\ell=0$ as well.

Suppose $\nu$ has a polynomial dependence in   $\omega$ and $\lim_{\omega\to 0}\nu(\omega)=\mb$.
Recall that ${\bf l}\in \C^k$ denotes the vector $(1,\ldots,1)$.
\begin{propn}
The function
$\sum_{j=1}^k \mu_j^\ell\:\tilde C_j\bigl(\nu(\omega),\mu+\frac{t}{\omega},\omega\bigr)$
is a polynomial in all its arguments.
Its specialization $\omega=0$  is a polynomial
\be
\label{Tr_explicit}
\vartheta_\ell(\mb,\mu,t)&=&\sum_{s=1}^{k} t^{s-1} \sum_{1\leq j_1<\ldots <j_s \leq k}m_{j_1}\ldots m_{j_s}
\sum_{{\bf d}\in \{\ell+1-s:s\}}\mu^{d_1}_{j_1}\ldots \mu^{d_s}_{j_s}.
\ee
\end{propn}
\begin{proof}
It seen from (\ref{ConnesI}) that the functions $\tilde C_j\bigl(\nu(\omega),\mu+\frac{t}{\omega},\omega\bigr)$ are
regular at $\omega=0$; this proves the first assertion.
Consider the specialization at $\omega=0$:
\be
\label{ConnesI1}
C_j(\mb,\mu,t)=\tilde C_j\bigl(\nu(\omega),\mu+\frac{t}{\omega},\omega\bigr)|_{\omega=0}.
\ee
We assume $\mu_i\not=0$ for all $i$ and put $\mb/\mu=(\mb_1/\mu_1,\ldots,\mb_k/\mu_k)$.
Now the statement follows from formulas (\ref{tr_explicit}) and (\ref{thetas}) if one observes  that
$
C_j(\mb,\mu,t)=
\tilde C_j\bigl(\mb/\mu,\mu,t\bigr)
$ from (\ref{ConnesI}).
\end{proof}
\noindent
Specialization at
$\nu=\hat \mb$, $\omega=1-q^{-2}$ yields polynomials $\vartheta_\ell(\hat \mb, \mu, q^{-2},t)$
participating in the two-parameter $\U_q(\g)$-equivariant quantization of orbits, \cite{DM2}:
\be
\label{qtConnesI}
\vartheta_\ell(\hat \mb, \mu, q^{-2},t)&=&
\sum_{j=1}^k \mu_j^\ell\:C_j(\hat \mb, \mu, q^{-2},t),
\ee
where $C_j(\hat \mb, \mu, q^{-2},t)=\tilde C_j\bigl(\hat \mb,\mu+\frac{t}{\omega},\omega\bigr)|_{\omega=1-q^{-2}}$.
The coefficients $C_j(\hat\mb,\mu,q^{-2},t)$ were shown above to give the non-commutative
Connes index for basic homogeneous vector bundles over the two-parameter quantum orbits.

\end{document}